\documentclass[11pt]{article}

\usepackage{hyperref}
\usepackage{harvard}
\usepackage{textcomp}
\usepackage{amsmath}
\usepackage{amsfonts}
\usepackage{epsfig}
\usepackage{color}

\usepackage{amsthm}
\usepackage{url}
\usepackage{graphicx}
\usepackage{multicol}
\usepackage[table]{xcolor}
\usepackage{enumerate}
\usepackage{graphicx}
\usepackage{subfig}
\usepackage[left=3cm, right=3cm]{geometry}
\usepackage{multirow}

\newcounter{gro}
\newenvironment{grouping}
{\begin{list}{Group \Roman{gro} :}
{\usecounter{gro}}}
{\end{list}}

\renewcommand{\baselinestretch} {1.3}

\makeatletter 
\def\singlespace{\def\baselinestretch{1}\@normalsize}

\title{{\sc Short-Term Load Forecasting: The Similar Shape Functional Time Series Predictor}}

\author{ Efstathios Paparoditis and  Theofanis Sapatinas\vspace{0.2cm}\\
Department of Mathematics and Statistics,  University of Cyprus, \\P.O. Box 20537, CY 1678 Nicosia, Cyprus.\\
 Email: paparoditis@ucy.ac.cy ~ Email: fanis@ucy.ac.cy}

\date{}

\newcommand{\reals}{\ensuremath{{\mathbb R}}}
\newcommand{\RR}{\reals}

\newcommand{\NN}{\ensuremath{{\mathbb N}}}
\newcommand{\EE}{\ensuremath{{\mathbb E}}}

\newcommand{\II}{\ensuremath{{\mathbb I}}}

\newtheorem{theorem}{Theorem}[section]
\newtheorem{lemma}{Lemma}[section]
\newtheorem{remark}{Remark}[section]

\def\endproof{\mbox{\ $\Box$}}

\begin{document}

\DeclareGraphicsExtensions{.pdf, .jpg}

\maketitle

\begin{abstract}
We introduce a novel functional time series methodology for short-term load forecasting. The prediction is performed by means of a weighted average of past daily load segments, the shape of which  is similar to the expected shape of the load segment to be predicted. The past load segments are identified from the available history of the observed load segments by means of their  closeness to a so-called reference load segment, the later being selected in a manner that captures the expected qualitative and quantitative characteristics of the load segment  to be  predicted. Weak consistency of the suggested functional similar shape predictor is established. As an illustration, we apply the suggested functional time series forecasting methodology to historical daily  load data in Cyprus and  compare its performance to that of a recently proposed alternative functional time series methodology for short-term load forecasting.
\medskip

\vspace{0.2cm}

{\bf AMS(2010) Subject Classification:} Primary  60G25, 62M20; Secondary 62P30

\end{abstract}


\section{Introduction}
\label{sec:Intro}

Load forecasting is an integrable process in the design of power systems  faced by electricity authorities worldwide, involving accurate predictions of electric load over different time periods in the furure. It can be broadly classified as short-term, medium-term and long-term forecasting, in terms of planning time horizons. Although different planning time horizons for these categories seem to exist in the literature, we have adopted herein the load forecasting classification scheme of Srinivasan and Lee (1995), viz., up to 1 day for short-term load forecasting (STLF), more than 1 day up to 1 year for medium-term load forecasting (MTLF), and more than 1 year up to 10 years for long-term load forecasting (LTLF). See also Alfares \& Nazeeruddin (2002).

Competition, the need of saving raw materials that are used in the production of electrical energy, the reduction of emissions and the avoidance of money wasting in general, are the main reasons that force electricity authorities worldwide to proceed to a better planning for the production of electricity. The main characteristics of the programming are the quantity of electrical energy that needs to be produced and the type of machine that is going to be used. This can be achieved by requiring accurate STLF, an important category of load forecasting, that plays a major role in real-time control and security functions for designing larger power systems, see, e.g., Srinivasan (1998) and Laurent {\em et al.} (2007).
Additionally, the particular characteristics of the electricity production process are essential for optimal planning of daily power generation.  In connection with the fact that the produced electricity that is not consumed instantly is lost (electricity cannot be stored), it becomes obvious that the right planning of daily power generation must have as a result the avoidance of emergency situations as well as of producing much greater quantities of electricity than the ones needed. It must also offer the capability to electricity authorities worldwide to use as far as possible the low functionality cost machines for covering the electrical energy needed. As a result, it becomes easily understood that an accurate STLF, like the prediction of the consumption of electrical energy of the next day, is an important tool for good power planning.  Over the years, a large number of methodologies have been developed to perform STLF. These methodologies are mainly emerged from two different paradigms, namely classical  statistical techniques and computational intelligent techniques. In particular, classical statistical techniques include, among others, regression models, ARIMA time series models, Kalman filtering models, etc., see, e.g., Hyde \& Hodnett (1997), Sargunaraj, Gupta \& Devi (1997) and Huang \& Shih (2003), while computational intelligent techniques include, among others, artificial neural networks, expert systems, etc., see, e.g., Rahman \&d Hazim (1993) and AlFuhaid, El-Sayed \& Mahmoud (1997). For extensive reviews on classical statistical techniques and computational intelligent techniques for STLF, we refer to, e.g.,  Hippert, Pedreira \& Souza (2001), Alfares \& Nazeeruddin (2002), Kyriakides \& Polycarpou (2007) and Taylor \& McSharry (2008).

Here, we focus on statistical techniques for STLF. It should be stressed that STLF is commonly considered as a difficult task because the daily load demand is influenced by many factors, viz., weather conditions, holidays, weekdays, weekends,  economic conditions, and, last but not least, idiosyncratic and social habits of individuals. Moreover, daily load demand is commonly recorded at a finite number of equidistance time points, viz., every half of each hour or every quarter of each hour. Thus, in order to forecast the load demand of next day, one has to predict the load demand at forty-eight or ninety-six, respectively, time points. From a statistical point of view, it is convenient to think of the daily load demand recorded at these forty-eight or ninety-six time points as a {\em segment} and to perform load prediction for the whole segment of time points rather than forecasting the load demand at each one of these time points separately. This implies that we adopt the {\em functional}\, time series framework in our approach, see, e.g., Bosq (2000, Chapter 9) and Ferraty \& Vieu (2006, Chapter 12).

Based on the previous discussion, we introduce below a novel functional time series methodology for STLF. The prediction is performed by means of a weighted average of past daily load segments, the shape of which  is similar to the expected shape of the load segment to be predicted. The past load segments are identified from the available history of the observed load segments by means of their  closeness to a so-called reference load segment, the later being selected in a manner that captures the expected qualitative and quantitative characteristics of the load segment  to be  predicted.

The paper is organized as follows. In Section \ref{sec:metho}, we provide some methodological background and review alternative functional time series methodologies that can be applied for STLF. We then describe the suggested functional time series forecasting methodology and establish weak consistency of the resulting functional predictor. As an application, we illustrate in Section \ref{sec:EAC}, how the proposed functional predictor performs by applying it 
to the historical daily electricity load data in Cyprus and compare its performance to that of a recently proposed alternative functional time series methodology for STLF. Some concluding remarks are provided in Section \ref{sec:con}.

\section{Functional Time Series Forecasting}
\label{sec:metho}

Putting the above discussion in a statistical methodological context, one seeks information on the evolution of a (real-valued) continuous-time stochastic process $X = (X(t); \; t \in \reals)$
in the future. Given a trajectory (curve) of $X$ observed on the interval $[0, T]$, one
would like to predict the behavior of $X$ on the entire interval
$[T,T+\delta]$, where $\delta > 0$, rather than at specific time-points. An
appropriate approach to this problem is to divide the interval $[0, T]$ into
subintervals $[l\delta,(l+1)\delta]$, $l=0,1,\ldots,k-1$ with $k=T/\delta$, and
to consider the (function-valued) discrete-time stochastic process ${\cal S} = ({\cal S}_n; \; n \in \NN)$, where
$\NN=\{1,2,\ldots\}$, defined by
\begin{equation}
{\cal S}_n(t) = X(t+ (n-1)\delta), \quad n \in \NN, \quad \forall \;t \in [0,
\delta). \label{eqn:zprocess}
\end{equation}
For the specific STLF application in mind, where the aim is one-day ahead prediction, the segmentation parameter $\delta$ corresponds to the daily electricity demand. In practice, the electricity demand is recorded at a finite number of equidistance time points within each day, say  $t_1,t_2,\ldots,t_P$, for instance, every half of each hour (viz., $P=48$) or every quarter of each hour (viz., $P=96$). Letting $S_n(t_i)$ to be the observation at time point $t_i$, $i=1,2,\ldots,P$, within curve ${\cal S}_n$,  $n \in \NN$, we denote by
$$
S_n =[S_n(t_1), S_n(t_2), \ldots, S_n(t_P)], \quad n \in \NN,
$$
the segment of the total number of observations of the $n$-th curve ${\cal S}_n$,  $n \in \NN$. 

Therefore, given a `sample' $S_1,S_2,\ldots,S_L$ of segments, our aim is then to predict the whole next segment $S_{L+1}$, viz., to predict
$$
S_{L+1} =[S_{L+1}(t_1), S_{L+2}(t_2), \ldots, S_{L+1}(t_P)].
$$

\subsection{Existing Approaches}
\label{subsec:exist}

In the recent statistical literature, practically all investigations to date for the aforementioned
functional time series prediction problem are for the case where one assumes that the underlying stochastic process is driven by a Hilbert-valued, first-order, autoregressive processes, which implies that the best predictor, ${\hat {\cal S}}_{L+1}$, of curve
${\cal S}_{L+1}$ given its past history (the `sample' of curves) ${\cal S}_1,{\cal S}_{2},\ldots,{\cal S}_L$, is the conditional mean of ${\cal S}_{L+1}$ given ${\cal S}_L$, see,
e.g., Bosq (1991), Bosq (2000, Chapter 3). In practice, however, an appropriate version of ${\hat {\cal S}}_{L+1}$ is obtained using some regularization on the predictor ${\hat S}_{L+1}$ of segment $S_{L+1}$. In particular, projection, spline and wavelet based regularization techniques have been proposed and consistency results of the resulting predictions have been obtained, see, e.g., Bosq (2000, Chapter 4), Besse \& Cardot (1996), Pumo (1998), Besse, Cardot \& Stephenson (2000) and Antoniadis \& Sapatinas
(2003). For some recent variants of this autoregressive structure, see, e.g., Marion \& Pumo (2004) and Mas \&  Pumo (2007).

An alternative approach to this prediction problem,  was proposed by Antoniadis, Paparoditis \& Sapatinas (2006). They developed a predictor via functional kernel nonparametric regression estimation
techniques, using a conditioning idea. In particular, prediction of segment $S_{L+1}$ was obtained by kernel smoothing, conditioning on the last observed segment $S_L$, while the resulting predictor was expressed as a weighted average of the past segments, placing more weight on those segments the preceding of which is similar to the present one. This functional time series forecasting methodology is rooted in the ability to find `similar' segments. Considering that segments can be sampled values of quite irregular curves, similarity matching was based on a distance metric on the discrete wavelet coefficients
of a suitable wavelet decomposition of the available segments.   For a similar approach, see, e.g., Ferraty,
Goia \& Vieu (2002) and Ferraty \& Vieu (2006, Chapter 11).

The common implicit assumption in developing the above functional time series forecasting methodologies is that all the available information for predicting segment $S_{L+1}$ is essentially contained in the last observed segment, viz., segment $S_L$. However, it is more appropriate to assume that the profile of the daily electricity load demand depends, in a complicated and unknown way, on a set of quantitative variables, such as,  daily temperature, daily humidity, daily wind speed, etc.,  as well as on a set of qualitative variables, such as, weekdays, weekends,  holidays, seasonal characteristics, etc., of the {\em day to be predicted}. Apparently, this information is not necessarily contained in the behavior of the last (observed) segment $S_L$. Thus, functional time series forecasting approaches that are based on conditioning ideas on the past behavior of the observed segments can ignore important information concerning the segment to be predicted, which is not necessarily contained in the past observed segments. The numerical results we present in Section \ref{sec:EAC} clearly show that this information is proved  valuable for accurate daily load demand forecasting.

Recently, Ferraty {\em et al.} (2011) considered a kernel regression estimator when both the response and the explanatory variables are functional. Thinking of the explanatory variables being quantitative, such as, daily temperature, daily humidity, daily wind speed, etc., the resulting functional kernel regression approach could be used for STLF. However, this functional regression-based approach does not take appropriately into account a number of  specific factors, the behavior of which turns out to be important for STLF. For instance, the same behavior of the daily temperature as an explanatory variable could lead to a different response, viz., daily load demand,  depending on the seasonal characteristics and on other factors, such as, weekdays, weekends,  holidays,  etc. This suggests that in order to perform accurate STLF, one needs to appropriately take into account not only the behavior of some quantitative variables but also on qualitative characteristics of the segment to be predicted, such as, weekdays, weekends,  holidays, and seasonal factors, that jointly affect the daily load demand behavior.

The above results in typical curves of daily load demand behavior depending on a number of quantitative and qualitative variables. Thus, in order to perform STLF one has to: (i) identify the appropriate curve of next day's load demand based on next day's behavior of these variables, and (ii) find in the entire time series history those curves that are similar to the identified one. With this in mind, we introduce below an alternative functional time series forecasting methodology for daily electricity load demand which performs the prediction by means of a weighted average of past segments, the shape of which  is similar to the expected shape of the segment to be predicted. The past segments are identified from the entire history of the observed segments  by means of their  closeness to a so-called reference segment, the later being selected in a manner that captures the expected qualitative and quantitative characteristics of the segment to be  predicted. Since the proposed STLF methodology is looking at the entire past for `shapes' that are similar to the expected `shape' of the day to be predicted, the resulting predictor is called the {\em similar shape functional time series predictor} (SSP).

\subsection{The Similar Shape Functional Time Series Predictor}
\label{subsec:Methodology}
Based on the previous discussion, it is not unrealistic to assume that daily load demand depends on several quantitative factors, such as, daily temperature, daily humidity,  daily wind speed, etc., and qualitative variables, such as, weekdays, weekends, holidays, season characteristics, etc.  Suppose that these quantitative and qualitative characteristics result in a typical daily load demand shape, so that each curve ${\cal S}_n$ can be expressed as
\begin{equation}
\label{eq:funcmod}
{\cal S}_n(t) = \sum_{m=1}^{M} f_m({\cal T}_n(t), {\cal H}_n(t), {\cal W}_n(t),\ldots)\, \II(g_n \in {\cal G}_m) + \varepsilon_n(t), \quad n \in \NN, \quad t \in [0, \delta),
\end{equation}
where
\begin{itemize}
\item  $ M$ is an {\em unknown}  number of typical daily load demand shapes described by the {\em unknown} functions  $f_m$ that correspond to different values of ${\cal G}_m$,  $m=1,2, \ldots,M$.

\item ${\cal G}_m=({\cal G}_{1,m}, {\cal G}_{2,m},\ldots)$ denotes a group of {\em deterministic} variables, denoting qualitative characteristics, e.g.,
${\cal G}_{1,m}$ refers to  {\em grouping of days} (weekdays, weekends, holidays), ${\cal G}_{2,m}$ to the {\em season characteristics}, etc. 

(Here, and in what follows, $\II(A)$ denotes the indicator function of the set $A$.)

\begin{remark}{\rm Notice that, in order to identify the seasonal characteristics, we do not only take into account the `global seasonality', as this is described by autumn, winter, spring and summer, but also the `local seasonality',  which refers to the weather conditions during the very recent past and which affect the behavior of the daily load demand. This is important since, for instance, a period of warm days  in the winter causes a different behavior of the daily load demand compared to the one caused by a period of similar warm days in the spring or in the summer.  This local seasonality aspects are also one of the reasons why regression based approaches are not very appropriate in this context. For instance, the same value of the daily temperature (the explanatory variable) may cause a different daily load demand (the response variable), depending on these local seasonal characteristics.}
\end{remark} 

\item $ g_n=(g_{1,n},g_{2,n},\ldots)$ denotes the qualitative characteristics of curve ${\cal S}_n$, $n \in \NN$, e.g., $g_{1,n}$ refers to the  {\em group membership} (weekdays, weekends, holidays), $g_{2,n}$ refers to the {\em season membership} (autumn, winter, spring, summer), etc.

\item ${\cal T}_n(t), {\cal H}_n(t), {\cal W}_n(t),\ldots $ are exogenous random variables, denoting quantitative characteristics, e.g.,
${\cal T}_n(t)$ is the {\em daily temperature} curve, ${\cal H}_n(t)$ is the {\em daily humidity} curve, ${\cal W}_n(t)$ is the {\em daily wind speed} curve, etc. 

Notice that these curves are function-valued random variables, usually called {\em functional} random variables in the statistical literature. In other words, they are random variables defined on a common probability space with values in an infinite-dimensional metric space.

\item $(\varepsilon_n(t),\, n \in \NN)$ is a $C([0,\delta))$-valued
strong Gaussian white noise, viz., a sequence of independent and
identically distributed (i.i.d) $C([0,\delta))$-valued Gaussian
random variables, where $C([0,\delta))$ is the space of continuous functions defined on the interval $[0,\delta)$. Furthermore, it is assumed that $\EE(\varepsilon_n) =0$ and $\EE(\|\varepsilon_n\|^2) <
\infty $. (Note that, in this case, the errors $\epsilon_n(t_i)$,
$i=1,2,\ldots,P$, forms a sequence of i.i.d. Gaussian random
variables with zero mean and finite variance.) It is also assumed that $\varepsilon_n(t)$ is
independent of ${\cal T}_s(t), {\cal H}_s(t), {\cal W}_s(t),\ldots $ for all $ s\leq n$, $n \in \NN$.

\end{itemize}

In what follows, for simplicity and data availability issues, we restrict our attention to the daily temperature curve ${\cal T}_n$, $n \in \NN$, which is one of the main exogenous functional random variables affecting daily load demand, that is, we consider the following model
\begin{equation}
\label{eq:funcmodNew}
{\cal S}_n(t) = \sum_{m=1}^{M} f_m({\cal T}_n(t))\, \II(g_n \in {\cal G}_m) + \varepsilon_n(t), \quad n \in \NN, \quad t \in [0, \delta).
\end{equation}
(However, we note that the suggested methodology can be straightforwardly modified to any available number of exogenous functional random variables.) Let $T_{L+1}(t_i)$
to be the observation at time point $t_i$, $i=1,2,\ldots,P$, within curve ${\cal T}_{L+1}$, viz., we denote by
$$
T_{L+1} =[T_{L+1}(t_1), T_{L+1}(t_2), \ldots, T_{L+1}(t_P)]
$$
the segment of the total number of observations of the $(L+1)$-th curve ${\cal T}_{L+1}$. 

Given model (\ref{eq:funcmodNew}) and based on the `sample' $S_1,S_2,\ldots,S_L$ of segments, the suggested prediction procedure consists of the following two main steps: (a) identify among the $M$ typical daily load demand shapes the one which is appropriate for the prediction of segment $S_{L+1}$, depending on the quantitative and qualitative characteristics of the segment to be predicted, and (b) obtain the prediction as a weighted average of all past segments by giving more weight to the segments the shape of which is more similar to the identified shape (obtained in (a)) of the segment to be predicted.  

\medskip
The following algorithm describes in more detail how to construct the suggested SSP, $\hat{S}_{L+1}$, of segment $S_{L+1}$. In particular, Steps 1 and 2 sort out the problem of identifying the appropriate shape of the segment to be predicted described in (a) while Step 3 refers to the calculation of the predictor described in  (b).

\begin{enumerate}

\item[Step 1.] Specify the {\em group membership} ${\cal G}_{1,m}$ and the {\em local  season membership} ${\cal G}_{2,m}$, where segment $S_{L+1}$ belongs to. While specifying  the group membership
 is easily done based on the particular day to be predicted (weekdays, weekends, holidays), the  specification of  the `local seasonality' is more difficult and rather arbitrary.  This is so, since, as explained before, local seasonality
  depends on the specific seasonal and weather characteristics of the very recent past of the time series and their stability. In our approach, local seasonality is essential and is taken into account by
 selecting a small number $n_L$ of past segments that are further considered for selecting what it is called the typical shape or {\em reference segment}. This is done in Step 2 of the algorithm.

\item[Step 2.] Determine a relevant load profile by specifying a so-called reference segment,
$S^{(Re)}$, as follows:

\begin{itemize}
\item Find among the last $n_L$ segments those belonging to same {\em group memberships} ${\cal G}_{1,m}$, $m \in \{1,2,\ldots,M\}$ as segment $S_{L+1}$. Let ${\cal C}_{L+1}$ be the set of the selected segments, viz.,
\begin{equation}
\label{eq:clF}
{\cal C}_{L+1}= \{S_l: \, g_{1,l} = {\cal G}_{1,m},\, 1 \leq l \leq n_L\}.
\end{equation}
Notice that the length of local seasonality and its stability is controlled by the parameter $n_L$ which determines how far in the past one goes to select a possible set of segments for specifying the reference segment $S^{(Re)}$.

\item Let $\widehat{T}_{L+1}$ be a predictor of segment $T_{L+1}$, which should not necessarily be available on the entire set of time points $\{t_1,t_2,\ldots,t_P\}$ (see below).

\item Let ${\cal D}$ be any of the (equivalent) {\em distances} in $\RR^P$. Then, the {\em reference segment} $S^{(Re)}$ is obtained as
\begin{equation}
\label{eq:rs}
S^{(Re)}(t_i) = \frac{1}{|C^{\ast}|} \sum_{S_l \in C^{\ast}}  S_l(t_i),\quad i=1,2,\ldots,P,
\end{equation}
where
$$
C^{\ast} = \{S_l \in {\cal C}_{L+1}:\, \II({\cal D}(T_l, \widehat{T}_{L+1}) \leq \delta)\}
$$
and
\begin{equation}
\label{eq:deltminF}
\delta  :=\delta(L)  \geq  {\min}_{\{l \in {\cal C}_{L+1}\}}\{{\cal D}(T_l, \widehat{T}_{L+1})\}.
\end{equation}
Notice that the reference segment $S^{(Re)}$ is obtained as a simple average of those segments $S_l \in{\cal C}_{L+1}$ belonging to the same group and having the same local seasonal
characteristics, the temperature segments $T_l$, $l=1,2,\ldots,n_L$,  of which are close enough to the predicted temperature $\widehat{T}_{L+1}$  of the segment $S_{L+1}$ to be predicted.
This `closeness' is controlled by the parameter $ \delta$. We also point out (and this is important for practical applications) that it is not necessary to have the prediction $\widehat{T}_{L+1}$ on the entire set of time points $\{t_1,t_2,\ldots,t_P\}$ to specify the set $C^{\ast}$. In fact, one can compare the temperature segments $T_l$ and $\widehat{T}_{L+1}$ using only the subset of time points  on which the predictions for the segment $T_{L+1}$ are available (or provided by other sources). 

\end{itemize}

\item[Step 3.]Finally, the SSP, $\hat{S}_{L+1}$, of the segment $S_{L+1}$, is given by 
\begin{equation}
\label{eq:fssp}
\hat{S}_{L+1}(t_i) = \sum_{r=1}^L w_r\,S_r(t_i), \quad i=1,2,\ldots,P,
\end{equation}
where the weights $w_r=w(S_r, S^{(Re)})$, $r=1,2,\ldots,L$, satisfy $w_r \geq 0$,  $r=1,2,\ldots,L$,
and $\sum_{r=1}^L w_r =1$. Following the nonparametric literature, the weights $w_r$,  $r=1,2,\ldots,L$, are chosen as
\begin{equation}
\label{eq:weightnewF}
w_r= \frac{K_{h}({\cal D}(S_r,S^{(Re)}))}{\sum_{l=1}^{L}K_{h}({\cal D}(S_l,S^{(Re)}))},\quad r=1,2,\ldots,L,
\end{equation}
with $K_h(\cdot)=h^{-1}K(\cdot/h)$ for some kernel function $K$, bandwidth $h_L$ and distance measure ${\cal D}$ between segments. It is assumed below that $K$ is a compactly supported bounded symmetric density such that $\int u^2 K(u)du <\infty$, and that $h=h(L) \rightarrow 0$ as $L \rightarrow \infty$.
\end{enumerate}

\begin{remark}
{\rm Notice that the SSP $\hat{S}_{L+1}$ is obtained as a weighted average of past segments, where more weight $w_r$, $r=1,2,\ldots,L$, is placed on the segment the shape of which is similar (in terms of the
particular distance ${\cal D}$ used) to the shape of the reference segment $S^{(Re)}$.  This clarifies  the differences between the suggested approach and several other  approaches proposed in the literature that  are based on conditioning ideas. In the case of conditioning on the last observed segment $S_L$ (see, e.g., Antoniadis, Paparoditis \& Sapatinas (2006)), the predictor is obtained as a weighted average of past segments,  where  the weight given to a segment depends on its closeness  to the conditioning segment $S_L$. For the suggested approach, the role of the conditioning segment is taken over by the reference segment $S^{(Re)}$ which comprehensively contains all relevant information regarding the shape of the daily load demand to be predicted. Thus, the selection of the reference segment $ S^{(Re)}$ is essential for the quality of the predictor obtained. This is stated more precisely in Lemma \ref{le:1}
which, in fact, shows that weak consistency of the suggested SSP is equivalent to weak consistency of the selected reference segment $S^{(Re)}$.}
\end{remark}
\medskip

Regarding the behavior of the smoothing parameters involved in the suggested prediction procedure, we assume the following.
\medskip

\noindent
{\bf Assumption 2.1}
\begin{itemize}
\item[(i)] $n_L \rightarrow \infty$ as $L \rightarrow \infty$ such that $|C^{\ast}| \rightarrow \infty $, where $C^{\ast}$ is defined in Step 2.

\item[(ii)] $\delta \rightarrow 0$ such that $|C^{\ast}| \delta \rightarrow \infty$.

\item[(iii)] $h \rightarrow 0$ such that $h L \rightarrow \infty$.

\end{itemize}
Assumption 2.1(iii) is standard for weak consistency in non-parametric kernel estimation.  Assumption 2.1(i) requires that the number of segments taken into account to calculate the reference segment $S^{(Re)}$, viz., the number of segments belonging to the set $C^{\ast}$, grows as the sample size increases. Assumption 2.1(ii) requires that the bandwidth $\delta$, used for obtaining the reference segment $S^{(Re)}$, goes to zero in such a way that the number of segments  effectively used in calculating $S^{(Re)}$, viz., $|C^{\ast}| \delta$,  increases to infinity, which is also as standard assumption for weak consistency in non-parametric kernel estimation.

\medskip

The following theorem establish the weak consistency of the suggested SSP.

\begin{theorem}
\label{th:1} Assume model (\ref{eq:funcmodNew}) and that Assumption 2.1 is satisfied. Let $\hat{S}_{L+1}$ be defined by (\ref{eq:fssp}) and assume that 
\begin{equation}
\label{discEq}
S_{L+1}=f_m(T_{L+1})+\varepsilon_{L+1},
\end{equation}
for some $m \in \{1,2,\ldots,M\}$, where
$$
f_m(T_{L+1})=[f_m(T_{L+1}(t_1)),f_m(T_{L+1}(t_2)),\ldots,f_m(T_{L+1}(t_P))]
$$
and
$$
\varepsilon_{L+1}=[\varepsilon_{L+1}(t_1),\varepsilon_{L+1}(t_2),\ldots,\varepsilon_{L+1}(t_P)].
$$
Then 
\begin{equation}
\label{eq:12}
{\cal D}(\hat{S}_{L+1},f_m(T_{L+1}) \stackrel{{\cal P}}{\rightarrow} 0,\;\; \text{as}\;\; L \rightarrow \infty,
\end{equation}
where $\stackrel{{\cal P}}{\rightarrow}$ denotes convergence in probability.
\end{theorem}

\section{Application: EAC Daily Load Data}
\label{sec:EAC}

\subsection{Description of the Data Set}
\label{subsec:DD}

Electricity Authority of Cyprus (EAC) is the organization that is responsible for the generation, transmission and distribution of electricity in Cyprus. The target of EAC is to provide Cypriots with high quality of safe and reliable services and activities at competitive prices. 
EAC uses two types of machines to produce electricity. The first type of machine is a {\em steam turbine} that uses dynamic pressure generated by expanding steam to turn the blades of a turbine. Almost all large non-hydro plants use this system. About $80\%$ of all electric power produced in the world is by use of steam turbine. The advantages of using such a type of machine are the high overall cogeneration efficiencies of up to $80\%$, the wide range of possible fuels including waste fuel and biomass, the production of high temperature/pressure steam and the established technology. On the contrary, we can mark the low electrical efficiencies, the slow start up times, the poor part load performance and especially the need for expensive high-pressure boilers and other equipment.
The second type of machine is a {\em diesel engine} which uses an electrical generator. Diesel generating sets are used in places without connection to the power grid, as emergency power-supply if the grid fails. Of course, they are widely used not only for emergency power but also many of them have a secondary function of feeding power to utility grids either during peak periods or during periods with a shortage of large power generators. Although, we may say that the cost of their functionality is forbidding.

The best planning for EAC is to avoid these emergency situations, use machines which have low functionality cost and high electrical efficiencies for electricity generation and distribution, which would cover the needs of the whole island. The most important characteristic of the right plan should be the fact that the quantity of electricity produced must not be greater than one needs since the additional electricity produced can not be stored and is lost. An extremely useful tool for the right plan is the accurate prediction of the consumption of the electrical energy of the next day.

\begin{figure}[h!]
\begin{center}
\includegraphics[trim=35mm 100mm 35mm 100mm,width=0.75\textwidth]{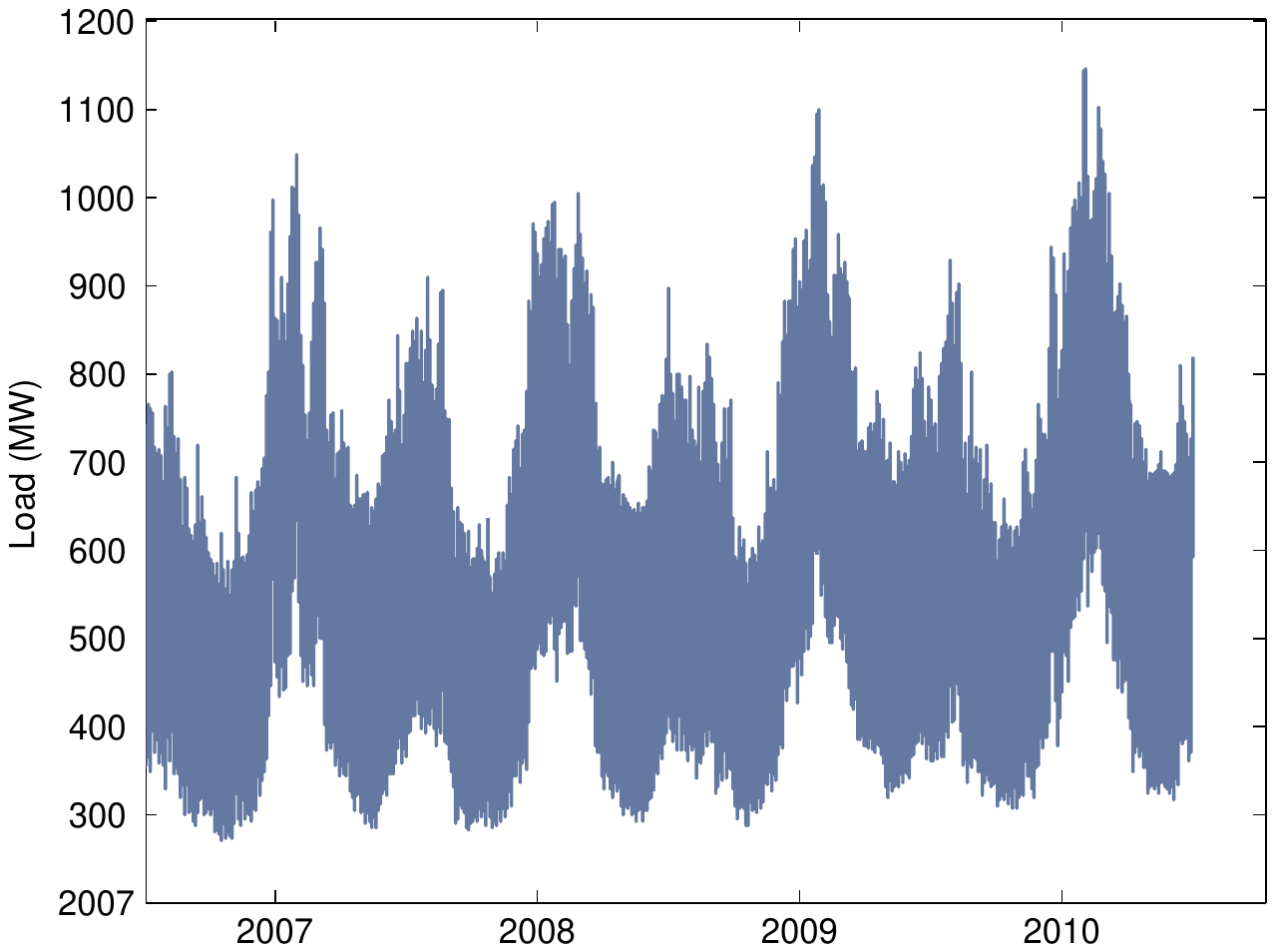}
\caption{Electrical power consumption in Cyprus between 1st January 2007 and 1st January 2010, recored every fifteen minutes.}
\label{fig:Th0}
\end{center}
\end{figure}

Below, we apply the suggested SSP methodology proposed in Subsection \ref{subsec:Methodology}
to a set of daily load data that were provided by EAC, concerning the electrical energy
consumption, in megawatts (MW), per fifteen minute intervals, viz., $P=96$, for the period from 01/01/2007 to 31/12/2010. The data set is displayed in Figure~\ref{fig:Th0}. From this figure, a slightly upward trend can be observed along with a strong periodic component within each year. It is also evident that the electrical energy  consumption slightly increases every year and during the summer months attains its maximum. Since the goal is to predict the daily shape of electrical energy consumption, and not the overall trend, the EAC rescales daily curves by dividing them by their maximum value. This leads to daily shape curves that vary between zero and one. The aim is then to produce accurate predictions of the shape of the rescaled consumption of next day's electrical energy. We then transform the predictor to the original scale, by multiplying the resulting SSP by the maximum value of the electrical load of the day to be predicted, provided by EAC. In the next section, we demonstrate how the SSP methodology proposed in Subsection \ref{subsec:Methodology}, can be implemented to fulfill this aim.

\subsection{Implementation of the SSP}
\label{subsec:imp}

In our context, the curves ${{\cal S}_n}$, $n = 1,2,\ldots,L$, that are derived from equation (\ref{eqn:zprocess}),
coincide with the calendar days from the 1st of January 2007 up to the last day ${{\cal S}_L}$, from which observations are available. Hence, ${{\cal S}_{L + 1}}$ coincide with the day for which prediction is required.
Based on the `sample' ${S_1},{S_2},\ldots,{S_L}$ of segments, the goal is to specify the SSP
${{\hat S}_{L + 1}}$. To this end,  the following steps are taken:

\begin{enumerate}
\item The group memberships ${\cal G}_{1,m}$, $m \in \{1,2,\ldots,M\}$, are specified,
where the segment $S_{L+1}$ belongs to. Feedback from the EAC, have shown that an appropriate  grouping of days, with similar shape behaviour based on some national characteristics, is the following:
\begin {grouping}
\item Monday, Tuesday, Thursday, Friday
\item Wednesday
\item Saturday
\item Sunday
\end{grouping}

\item To determine the reference segment $S^{(Re)}$, as mentioned previously, we
restrict our attention only to the exogenous random variable ${{{\rm \mathcal T}_n}\left( t \right)}$, $n \in \mathbb N$, that denotes the daily temperature segment. We then find all days belonging to same group memberships ${\cal G}_{1,m}$, $m \in \{1,2,\ldots,M\}$, as the segment we want to predict. The parameter $n_L$, described in Subsection \ref{subsec:Methodology}, is set equal to ${n_L} = 14$ if the segment to be predicted corresponds to Monday, Tuesday, Thursday or Friday and is set equal to ${n_L} = 28$ if the segment to be predicted corresponds to Wednesday, Saturday or Sunday. This seems to be appropriate in order to have enough information to select the reference segment while at the same time retains local seasonality. Regarding the exogenous random variable ${{{\rm \mathcal T}_n}\left( t \right)},n \in \mathbb N$, we use a predictor $\widehat{T}_{L+1}$ of $T_{L+1}$. As mentioned earlier, it is not necessary to have the prediction $\widehat{T}_{L+1}$ on the entire set of time points $\{t_1,t_2,\ldots,t_P\}$, $P=96$, to specify the set $C^{\ast}$. In fact, one can compare the temperature segments $T_l$, $l=1,2,\ldots,n_L$, and $\widehat{T}_{L+1}$ using only a subset of time points  on which the predictions for the segment $T_{L+1}$ are available (or provided by other sources). More specifically, we have used predictions of next day's temperature at only four time points, that is, those corresponding to 08:00, 12:00, 16:00 and 20:00. The actual temperatures predictions of segment $T_{L+1}$ at these time points were provided to us by EAC.

\item To further simplify the selection of the reference segment $S^{(Re)}$,  we set in (\ref{eq:deltminF})
$$
\delta  = {\min}_{\{l \in {\cal C}_{L+1}\}}\{{\cal D}(T_l, \widehat{T}_{L+1})\},
$$
that is, the reference segment $S^{(Re)}$ is obtained as
\begin{equation}
\label{eq:new2}
{S^{(Re)}}= {\arg\min}_{\{S_l \in {{\cal C}_{L + 1}}\}} \{ {\mathcal D( {{T_l},{{\widehat T}_{L + 1}}})}\},
\end{equation}
where ${\cal C}_{L + 1}$ is given by (\ref{eq:clF}).


\item The SSP $\hat{S}_{L+1}$ is obtained by (\ref{eq:fssp}),
where the weights ${w_r}$, $r=1,2,\ldots,L$, are determined by (\ref{eq:weightnewF}) with the kernel function $K$ being the Gaussian kernel and the bandwidth $h$ being selected by the empirical risk of prediction  methodology suggested by Antoniadis, Paparoditis \& Sapatinas (2009). 
\end{enumerate}


Finally, the computational algorithm related to the above implementation as well as the overall numerical study presented in Subection \ref{subsec:resul} has been carried out in the {\tt Matlab 7.7.0} programming environment.

\subsection{Numerical Results}
\label{subsec:resul}

Based on the above implementation, we apply the SSP functional time series forecasting methodology to a randomly selected number of days within the year 2010, from the dataset displayed in Figure~\ref{fig:Th0}, and compare its performance to the recently established wavelet-kernel functional time series methodology (WKP) proposed by Antoniadis, Paparoditis \& Sapatinas (2006). We restrict our comparison to that forecasting methodology, since as it has been demonstrated by these authors in a number of simulated and real-data examples, in terms of functional time series forecasting, WKP outperforms many
well-established forecasting methods, like a wavelet regularization method, a smoothing spline method, the classical SARIMA model and the Holt-Winters forecasting procedure. We refer to Antoniadis, Paparoditis \& Sapatinas (2006) for more details. 

The quality of both SSP and WKP were measured by the 
{\em relative mean-absolute error} (RMAE) defined as
\begin{equation*}
\label{eqn:rmae}
\text{RMAE} = \frac{1}{96} \sum_{i=1}^{96} \frac{|\hat{S}_{L+1}(t_i) -
S_{L+1}(t_i)|}{S_{L+1}(t_i)}.
\end{equation*}
For both SSP and WKP, we also report the {\em Maximum Difference} (MaxDiff) and {\em Minimum Difference} (MinDiff) defined as
$$
{\rm MaxDiff}= \max_{i=1,2,\ldots,96}\{\hat{S}_{L+1}(t_i) - S_{L+1}(t_i)\} \quad \text{and} \quad
{\rm MinDiff}= \min_{i=1,2,\ldots,96}\{\hat{S}_{L+1}(t_i) - S_{L+1}(t_i)\}.
$$

\begin{figure}[h!]
\begin{center}
\includegraphics[trim=35mm 100mm 35mm 100mm,width=0.75\textwidth]{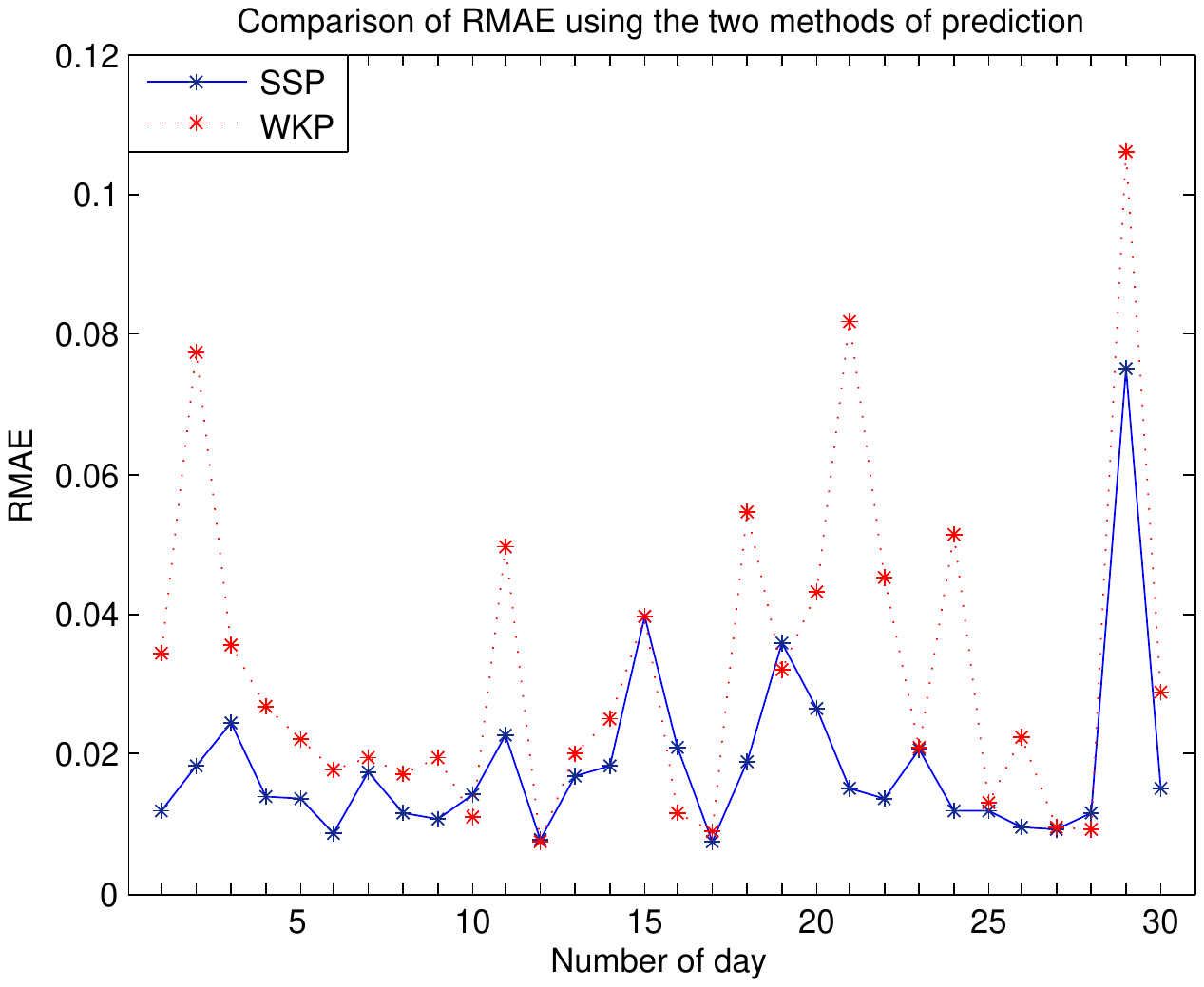}
\caption{RMAE for the Functional Similar Shape Time Series Predictor (SSP, solid line) and the Functional Wavelet-Kernel Time Series Predictor (WKP, dotted line) for the 30 randomly selected days within the year 2010.}
\label{fig:Th}
\end{center}
\end{figure}

\begin{figure}[h!]
\begin{center}
\includegraphics[trim=35mm 45mm 35mm 10mm,width=0.75\textwidth]{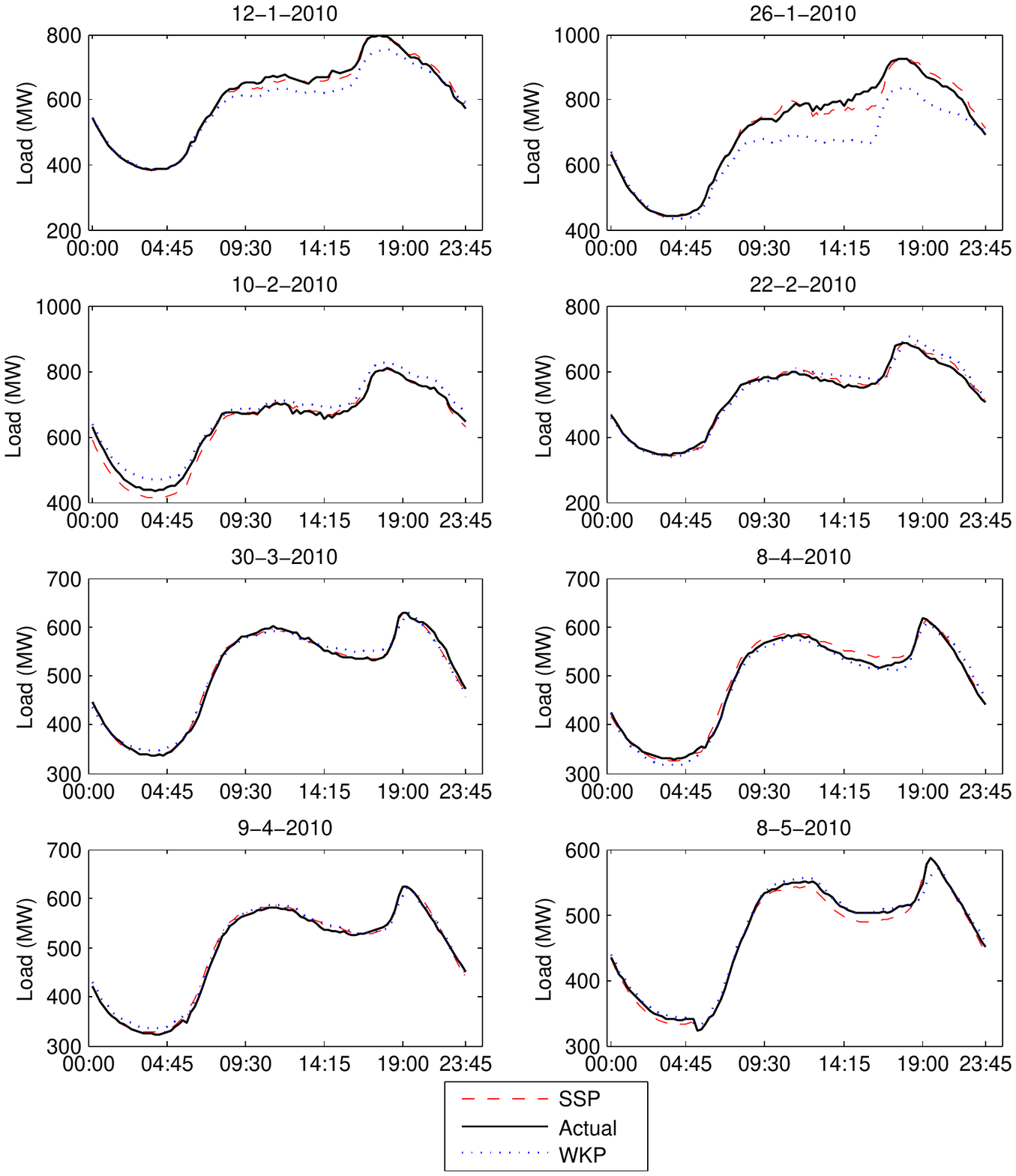}
\caption{Actual (solid line) and predicted load using the Functional Similar Shape Time Series Predictor (SSP, dashed line) and the Functional Wavelet-Kernel Time Series Predictor (WKP, dotted line) for the randomly selected days within the year 2010.}
\label{fig:Th1}
\end{center}
\end{figure}

\begin{figure}[h!]
\begin{center}
\includegraphics[trim=35mm 45mm 35mm 10mm,width=0.75\textwidth]{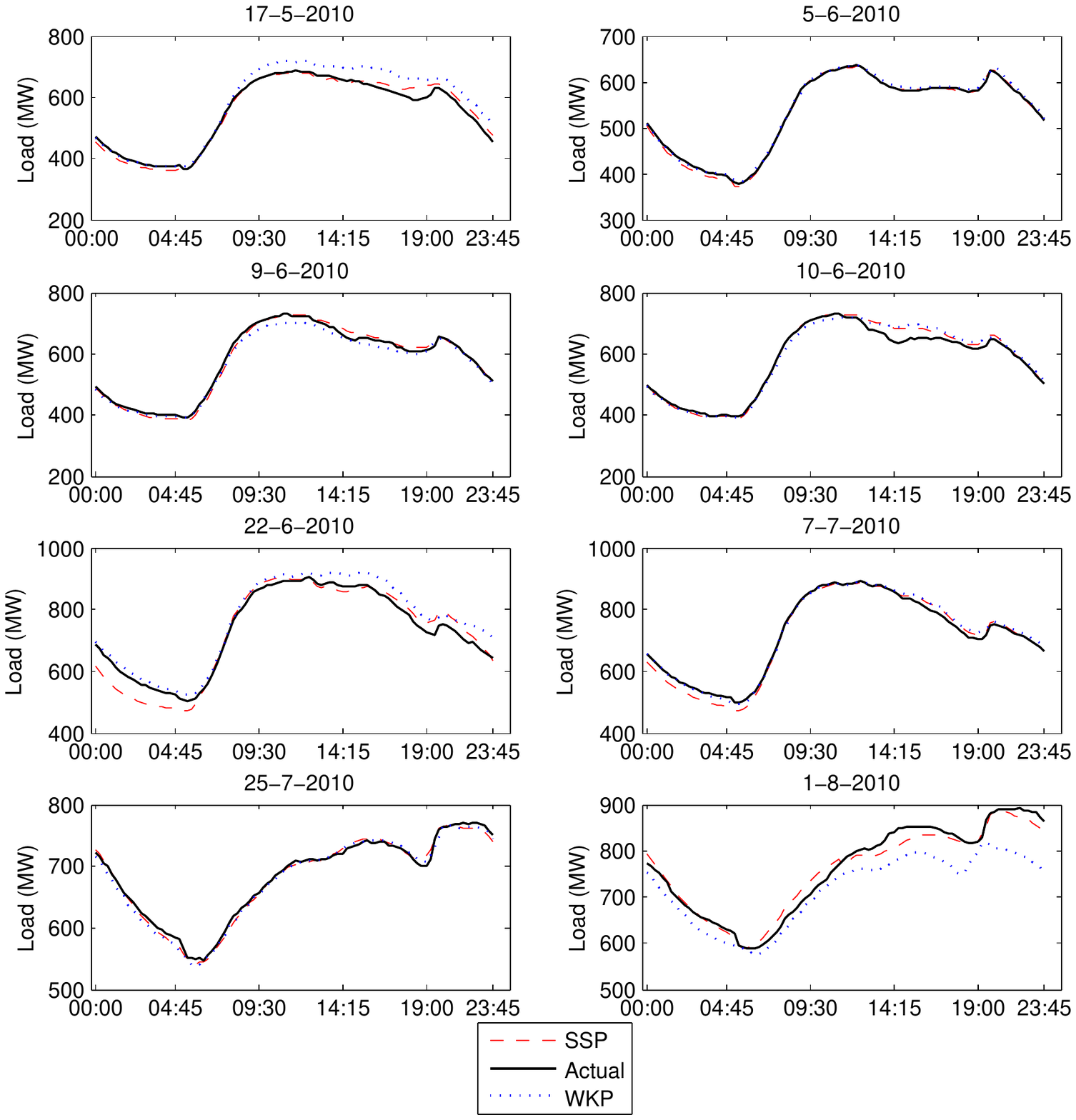}
\caption{Actual (solid line) and predicted load using the Functional Similar Shape Time Series Predictor (SSP, dashed line) and the Functional Wavelet-Kernel Time Series Predictor (WKP, dotted line)  for the randomly selected days within the year 2010.}
\label{fig:Th2}
\end{center}
\end{figure}

\begin{figure}[h!]
\begin{center}
\includegraphics[trim=35mm 45mm 35mm 10mm,width=0.75\textwidth]{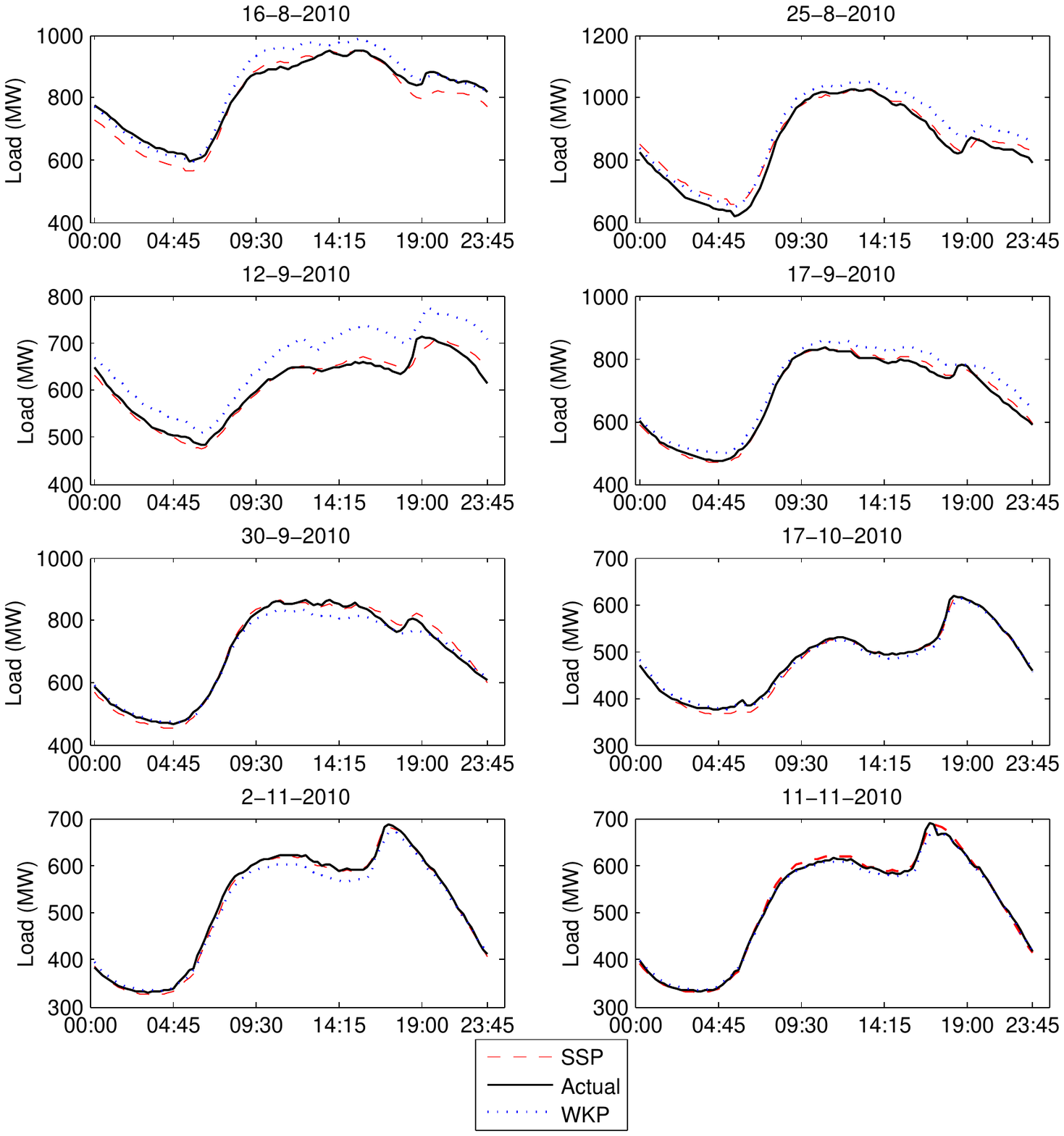}
\caption{Actual (solid line) and predicted load using the Functional Similar Shape Time Series Predictor (SSP, dashed line) and the Functional Wavelet-Kernel Time Series Predictor (WKP, dotted line) for the randomly selected days within the year 2010.}
\label{fig:Th3}
\end{center}
\end{figure}

\begin{figure}[h!]
\begin{center}
\includegraphics[trim=35mm 45mm 35mm 10mm,width=0.75\textwidth]{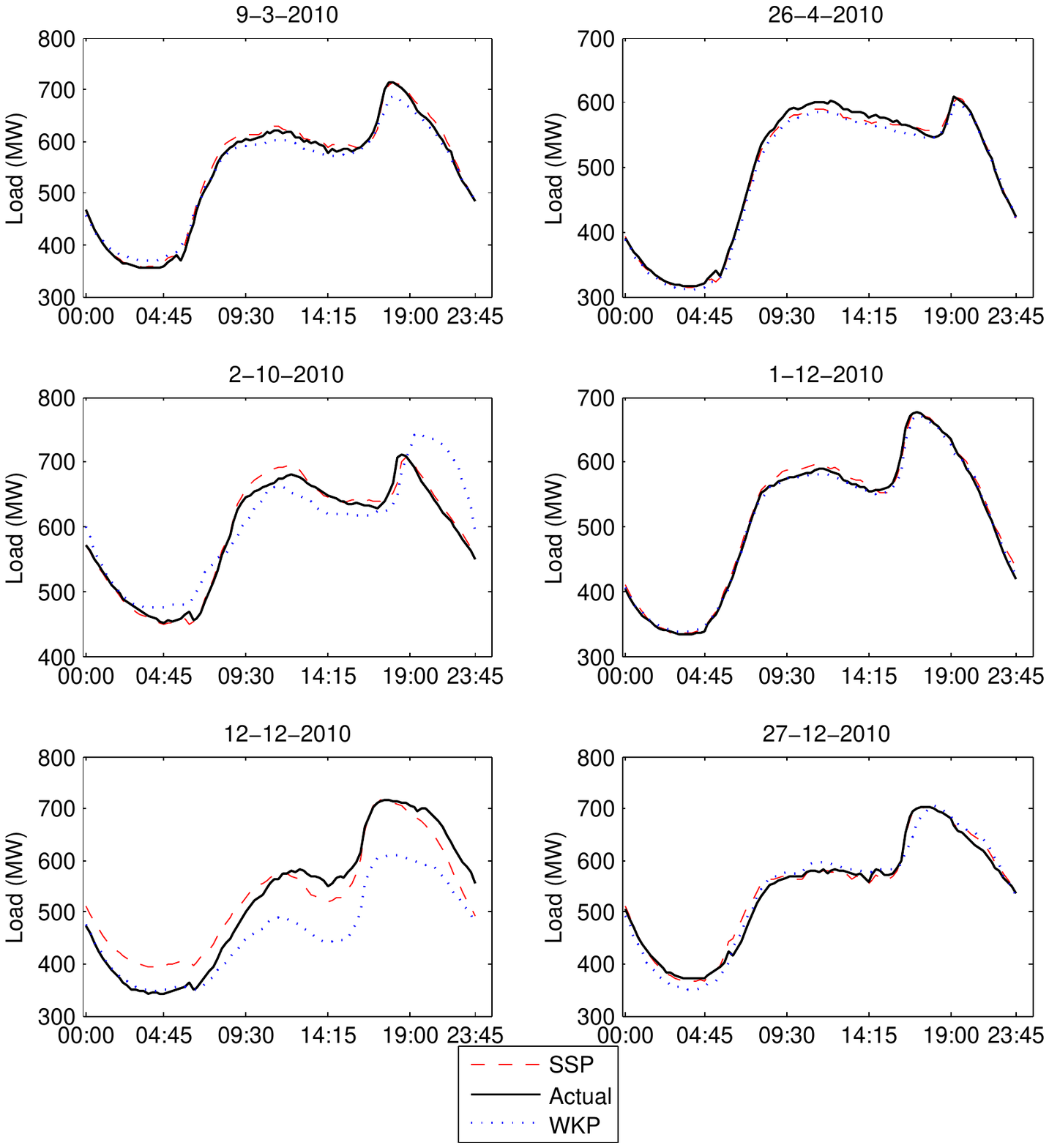}
\caption{Actual (solid line) and predicted load using the Functional Similar Shape Time Series Predictor (SSP, dashed line) and the Functional Wavelet-Kernel Time Series Predictor (WKP, dotted line) for the randomly selected days within the year 2010.}
\label{fig:Th4}
\end{center}
\end{figure}

It is evident from the analysis (see Figures \ref{fig:Th} --\ref{fig:Th4} and Table \ref{tab}) that the SSP 
clearly outperforms the WKP. It terms of RMAE, only in 5 out of the 30 randomly selected number of days the SSP performs slightly worse than the WKP. Furthermore,  in a large number of days, the RMAE of the WKP considerably exceeds the RMAE of the SSP, as it is clearly seen in Figure \ref{fig:Th}. Looking at each day separately (see Figures \ref{fig:Th1}--\ref{fig:Th4}), the SSP curves are quite close to the actual load curve while at the same time there are days where the WKP fails to appropriately capture even the overall behavior of the latter curve (see, e.g., 26 Jan 2010, 1 Aug 2010, 2 Oct 2010 and 12 Dec 2010).

\section{Conclusions}
\label{sec:con}

We introduced a novel functional time series methodology for short-term load forecasting, named the {\em functional similar shape time series predictor}. The predictor was performed by means of a weighted average of past daily load segments, the shape of which is similar to the expected shape of the load segment to be predicted.  To quantify this similarity, the notion of {\em reference segment} was introduced which captures 
the expected qualitative and quantitative characteristics of the load segment  to be  predicted. The functional similar shape predictor was theoretically justified by proving a weak consistency property. Furthermore, its usefulness for short-term load forecasting was demonstrated by applying it to historical daily  load data in Cyprus.  The numerical results obtained showed that the functional similar shape predictor works very satisfactory and outperforrms the functional wavelet-kernel time series predictor, a recently established alternative functional time series methodology for short-term load forecasting. 

We note that, although for simplicity and data availability issues, we restricted our attention to the daily temperature, which is one of the main exogenous functional random variables, the suggested functional time series methodology for short-term load forecasting can be  also modified to take into account other daily exogenous functional random variables, like humidity, wind speed, sunshine, etc., which might affect daily load demand. Although the above modification is straightforward from a theoretical point of view, its practical implementation depends on the availability of the required time series data.

\section{Appendix: Proof of Theorem \ref{th:1}}
\label{sec:app}

To prove Theorem \ref{th:1}, we first prove the following auxiliary result, showing that weak consistency of the suggested SSP $\hat{S}_{L+1}$, defined by (\ref{eq:fssp}), is equivalent to the weak consistency of the selected reference segment $S^{(Re)}$, defined by (\ref{eq:rs}).

\begin{lemma}
\label{le:1} Assume model (\ref{eq:funcmodNew}), that Assumption 2.1 $(iii)$ is satisfied, and that  (\ref{discEq}) holds true. Let $S^{(Re)}$ and $\hat{S}_{L+1}$ be defined by (\ref{eq:rs}) and (\ref{eq:fssp}), respectively. Then,  (\ref{eq:12}) holds true, if and only if
\begin{equation}
\label{eq:21}
{\cal D}(S^{(Re)},f_m(T_{L+1})) \stackrel{{\cal P}}{\rightarrow} 0,\;\; \text{as}\;\; L \rightarrow \infty.
\end{equation}
\end{lemma}

\noindent{\bf Proof.} We first prove the ``if'' part of the lemma. Assume that (\ref{eq:21}) holds true. Since
\begin{align*}
{\cal D}(\hat{S}_{L+1},f_m(T_{L+1})) & \leq
{\cal D}(\hat{S}_{L+1}, S^{(Re)}) + {\cal D}(S^{(Re)},f_m(T_{L+1})) \\ := A_L + B_L.
\end{align*}
Since 
$$
A_L \leq \sum_{r=1}^{L} w_r\, {\cal D}(S_{r}, S^{(Re)})=O_p(h) \rightarrow 0 
\quad \text{and} \quad B_L \stackrel{{\cal P}}{\rightarrow} 0,
\quad \text{as} \quad L \rightarrow \infty,
$$
the ``if'' part of the lemma follows. 

We now prove the ``only if'' part of the lemma. Assume that (\ref{eq:12}) holds true.
Note that, for each $t \in \{t_1,t_2,\ldots,t_P\}$,
$$
|\hat{S}_{L+1}(t) -S^{(Re)}(t)| \leq \sum_{r=1}^{L} w_r\, |S_{r}(t) -S^{(Re)}(t)| = O_p(h).
$$
Also,
\begin{align*}
{\cal D}(S^{(Re)},f_m(T_{L+1})) & \leq
{\cal D}(S^{(Re)}, \hat{S}_{L+1}) + {\cal D}(\hat{S}_{L+1},f_m(T_{L+1})) \\ := C_L + D_L.
\end{align*}
Since
$$
C_L = O_p(h) \rightarrow 0 \quad \text{and} \quad D_L \stackrel{{\cal P}}{\rightarrow} 0, \quad \text{as} \quad L \rightarrow \infty,
$$
the ``only if'' part of the lemma follows. \hfill \endproof

\medskip

We are now ready to prove Theorem \ref{th:1}.
In view of Lemma \ref{le:1}, it suffices to prove (\ref{eq:21}). From (\ref{eq:rs}), it is easily seen that $S^{(Re)}$ has the expression of a non-parametric estimator with a uniform kernel. Thus, in view of Assumption 2.1 $(i)$ and $(ii)$,  (\ref{eq:12}) holds true by standard weak consistency arguments for functional non-parametric kernel estimators, see, e.g., Bosq (1998, Chapter 3).
\hfill \endproof

\section*{Acknowledgements}
This project was partially supported by the Electricity Authority of Cyprus (EAC) and the Transmission System Operator, Cyprus (TSO). The authors would like to thank EAC and TSO for providing the electric load data and for many fruitful discussions on forecasting daily load demand in Cyprus.  The authors would also like to thank Miss Maria Frangeskou, an MSc student in Applied Statistics, in the Department of Mathematics and Statistics at the University of Cyprus, Cyprus, for her help with some of the numerical results.

\newpage

\bigskip

\newpage

\begin{table}
\begin{center}
\begin{tabular}{|l|lll|lll|}
\hline
\hline
     &  &SSP &  & & WKP& \\
Date & RMAE & MaxDiff & MinDiff & RMAE & MaxDiff & MinDiff \\
\hline
12 Jan 2010  (Wed) & 0.0120 & 21.735 & -21.350 & 0.0344& 13.010 & -62.318 \\
26 Jan 2010 (Wed) & 0.0183  & 31.771&-56.341 & 0.0774& 16.698&-168.246 \\
10 Feb 2010  (Thu)  & 0.0244&20.679&-40.814&0.0355& 48.139& -7.649\\
22 Feb 2010 (Tue) & 0.0139 & 19.484&-25.830&0.0269&41.484&-38.054\\
9 Mar 2010 (Wed) & 0.0136& 17.372&-14.837& 0.0220&27.278&-43.017\\
30 Mar 2010 (Wed) & 0.0085&17.398&-17.925& 0.0177& 21.790&-18.634\\
8 Apr 2010 (Fri) &0.0174&23.191&-12.960&0.0194&60.346&-29.092\\
9 Apr 2010 (Sat) & 0.0116 &17.667&-11.217& 0.0170& 19.972& -17.742\\
26 Apr 2010 (Tue) & 0.0107 & 11.446&-16.421&0.0195& 15.724&-54.782\\
8 May 2010 (Sun) & 0.0143&12.152&-17.117&0.0110&11.436& -29.488\\
17 May 2010 (Tue) & 0.0227& 42.119&-18.699& 0.0496&69.632&-7.392\\
5 Jun 2010 (Sun) & 0.0079& 8,.261&-12.915& 0.0075&12.039&-7.439\\
9 Jun 2010 (Thu) & 0.0168 & 23.480& -19.095&0.0201 &10.070&-32.029\\
10 Jun 2010 (Fri) & 0.0183 & 46.482& -9.387& 0.0251&50.860& -25.818\\
22 Jun 2010 (Wed) & 0.0397 & 46.887& -74.603 & 0.0398& 72.108& 0.0328\\
7 Jul 2010 (Thu) & 0.0209&19.830&-36.639 & 0.0117& 25.420& -16.504\\
25 Jul 2010  (Mon) & 0.0076& 16.212 & -12.496 & 0.0089& 8.709& -19.766\\
1 Aug 2010 (Mon) & 0.0188& 29.011&-37.558&0.0546&-83.500& -117.146\\
16 Aug 2010 (Tue) & 0.0359& 21.442&-79.744& 0.0320& 68.523&-16.523\\
25 Aug 2010 (Thu) & 0.0264 & 45.570&-31.062&0.0431&72.600&11.020\\
12 Sep 2010 (Mon) & 0.0151&30.734&-31.385&0.0817&96.777&19.450\\
17 Sep 2010 (Sat) & 0.0135& 22.667&-18.267&0.0452&66.735&-0.768\\
30 Sep 2010 (Fri) & 0.0207&35.214&-17.545&0.0208&18.353&-47.553\\
2 Oct 2010 (Sun) & 0.0119&20.838&-37.730&0.0514&105.903&-58.231\\
17 Oct 2010 (Mon) & 0.0120& 2.152&-17.166& 0.0131& 17.502&-21.257\\
2 Nov 2010 (Wed) & 0.0096& 9.516&-16.643& 0.0223&11.081&-28.035\\
11 Nov 2010 (Fri) & 0.0091& 19.189&-11.575&0.0095&11.235&-27.899\\
1 Dec 2010 (Thu) & 0.0117& 17.036&-11.819& 0.0093& 9.633& -24.258\\
12 Dec 2010 (Mon) & 0.0752& 54.466& -68.558&0.10617&8.019& -128.209\\
27 Dec 2010 (Tue) & 0.0151 & 40.785&-14.061& 0.0289& 25.129&-51.986\\
\hline
\hline
\end{tabular}
\end{center}
\caption{RMAE, MaxDiff and MinDiff for the Functional Similar Shape Time Series Predictor (SSP) and the Functional Wavelet-Kernel Time Series Predictor (WKP) for the 30 randomly selected days within the year 2010.}
\label{tab}
\end{table}
\newpage

\newpage

\newpage

\newpage

\begin{thebibliography}{}
\bibitem{}
Alfares, H.K. \& Nazeeruddin, M. (2002). Electric load forecasting: literature survey and classification of methods. {\em International Journal of Systems Science}, {\bf 33}, 23--34.

\bibitem{}
AlFuhaid, A.S., El-Sayed, M.A. \& Mahmoud, M.S. (1997). Cascaded artificial neural networks for short-term load forecasting. {\em IEEE Transactions on Power Systems}, {\bf 12}, 1524--1529.

\bibitem{}
Antoniadis, A., Paparoditis, E. \& Sapatinas, T. (2006). A
functional wavelet-kernel approach for time series prediction. {\em
Journal of the Royal Statistical Society}, Series B, {\bf 68},
837--857.

\bibitem{}
Antoniadis, A., Paparoditis, E. \& Sapatinas, T. (2009).
Bandwidth selection for functional time series prediction. {\em Statistics and Probability Letters},
{\bf 79}, 733--740.

\bibitem{}
Antoniadis, A. \& Sapatinas, T. (2003). Wavelet methods for continuous-time
prediction using Hilbert-valued autoregressive processes. {\em Journal of
Multivariate Analysis}, {\bf 87}, 133--158.

\bibitem{}
Besse, P.C. \& Cardot, H. (1996). Approximation spline de la pr\'{e}vision d'un
processus fonctionnel autor\'{e}gressif d'ordre 1. {\em Canadian Journal of Statistics},
{\bf 24}, 467--487.

\bibitem{}
Besse, P.C., Cardot, H. \& Stephenson, D.B. (2000). Autoregressive forecasting
of some functional climatic variations. {\em Scandinavian Journal of Statistics}, {\bf 27},
673--687.

\bibitem{}
Bosq, D. (1991). Modelization, nonparametric estimation and prediction for
continuous time processes. In {\em Nonparametric Functional Estimation and
Related Topics}, Ed. G. Roussas, pp. 509--529, Nato ASI Series C, Vol. {\bf
335}, Dortrecht: Kluwer Academic Publishers.

\bibitem{}
Bosq, D. (1998). {\em Nonparametric Statistics for Stochastic Processes}. Lecture Notes in
Statistics, 2nd Edition, Vol. {\bf 110}, New York: Springer-Verlag.

\bibitem{}
Bosq, D. (2000). {\em Linear Processes in Function Spaces}. Lecture Notes in
Statistics, Vol. {\bf 149}, New York: Springer-Verlag.


\bibitem{}
Ferraty, F., Goia, A. \& Vieu, P. (2002). Functional nonparametric
model for time series: a fractal approach for dimension reduction.
{\em Test}, {\bf 11}, 317--344.

\bibitem{}
Ferraty, F. \& Vieu, P. (2006). {\em Nonparametric Functional Data
Analysis}. New York: Springer-Verlag.

\bibitem{}
Ferraty, F., Laksaci, A., Tadj, A. \& Vieu, P. (2011). Kernel regression with functional response. {\em Electronic Journal of Statistics}, {\bf 5}, 159--171.

\bibitem{}
Hippert, H.S., Pedreira, C.E. \& Souza, R.C. (2001). Neural networks for short-term load forecasting: a review and evaluation. {\em IEEE Transactions on Power Systems}, {\bf 16}, 44--55.

\bibitem{}
Huang, S.-J. \& Shih, K.-R. (2003). Short-term load forecasting via ARMA model identification including non-Gaussian process considerations. {\em IEEE Transactions on Power Systems}, {\bf 18}, 673--679.

\bibitem{}
Hyde, O. \& Hodnett, P.F. (1997). An adaptable automated procedure for short-term electricity load forecasting. {\em IEEE Transactions on Power Systems}, {\bf 12}, 84--94.

\bibitem{}
Kyriakides, E. \& Polycarpou, M. (2007). Short term electric load forecasting: a tutorial. In {\em Trends in Neural Computing}, Studies in Computational Intelligence, Vol. {\bf 35}, (Eds., K. Chen and L. Wang), Chapter 16, pp. 391--418, Springer-Verlag: Berlin.

\bibitem{}
Laurent, P., Fock, E. Randrianarivony, R.N. \& Manicom-Ramsamy, J.-F. (2007). Bayesian neural network approach to short time load forecasting. {\em Energy Conversion and Management}, {\bf 49}, 1156--1166.

\bibitem{}
Marion, J.M. \& Pumo, B. (2004). Comparaison des mod\`eles ARH(1) et ARHD(1) sur
des donn\'ees physiologiques. {\em Annales de l'I.S.U.P.}, {\bf 48}, 29--38.

\bibitem{}
Mas, A. \& Pumo, B. (2007). The ARHD model. {\em Journal of Statistical Planning and Inference},
{\bf 137}, 538--553.

\bibitem{}
Pumo, B. (1998). Prediction of continuous time processes by $C[0,1]$-valued
autoregressive processes. {\em Statistical Inference for Stochastic Processes}, {\bf 3},
297--309.

\bibitem{}
Rahman, S. \& Hazim, O. (1993). A generalized knowledge-based short-term load-forecasting
technique. {\em IEEE Transactions on Power Systems}, {\bf 8}, 508--514.

\bibitem{}
Srinivasan, D. (1998). Evolving artificial neural networks for short-term load forecasting. {\em Neurocomputing}, {\bf 23}, 265--276.

\bibitem{}
Srinivasan, D. \& Lee, M.A. (1995). Survey of hybrid fuzzy neural approaches to electric load forecasting.
{\em Proceedings of the IEEE International Conference on Systems, Man and Cybernetics}, Part 5, Vancouver, BC, pp. 4004-4008.

\bibitem{}
Sargunaraj, S., Gupta., D.P.S. \& Devi, S. (1997). Short-term load forecasting for demand side management. {\em IEE Proceedings on Generation, Transmission and Distribution}, {\bf 144}, 68--74.

\bibitem{}
Taylor, J.W. \& McSharry, P.E. (2008). Short-term load forecasting methods: an evaluation based on European data. {\em IEEE Transactions on Power Systems}, {\bf 22}, 2213--2219.


\end{thebibliography}
\end{document}